\documentclass{article} 
\usepackage[usenames,dvipsnames]{color}
  \usepackage{amsmath} 
  \usepackage{amssymb} 
  \usepackage{bm}
  \usepackage{mathtools}
  \usepackage{amsthm}
 \usepackage{CJK}
 \usepackage{indentfirst}
\usepackage{tikz}
\usepackage{amsmath,amsfonts,amsthm,mathrsfs}
\numberwithin{equation}{section}
\newtheorem{theorem}{Theorem}
\newtheorem{lemma}{Lemma}

\newtheorem{definition}{Definition}

\usepackage{mathtools}
\usepackage{relsize}
\usepackage{multirow}
\def\NN{\mathbb N}
\def\ZZ{\mathbb Z}
\def\RR{\mathbb R}

\def\TT{\mathbb T}

\title{Generalization Analysis for Classification on Korobov Space}
\date{}
\author{Yuqing Liu \\
Department of Mathematics, HKU
 \\
Email: yuliu6@hku.hk\\
}

\begin{document}
\maketitle
\begin{abstract}
    \noindent
In this paper, the classification algorithm arising from Tikhonov regularization is discussed. The main intention is to derive learning rates for the excess misclassification error according to the convex $\eta$-norm loss function $\phi(v)=(1 - v)_{+}^{\eta}$, $\eta\geq1$. Following the argument, the estimation of error under Tsybakov noise conditions is studied. In addition, we propose the rate of $L_p$ approximation of functions from Korobov space $X^{2, p}([-1,1]^{d})$, $1\leq p \leq \infty$, by the shallow ReLU neural network. This result consists of a novel Fourier analysis approach and a probability argument.
\end{abstract}

\section{Introduction}

The challenge for misclassification problem in practice is that as the dimension grows large, the feature becomes into special forms. Therefore, a special structure contribution is required. The performance of classification of functions from Korobov space using shallow networks might be one of the possibilities to deal with the well.

A binary classification problem with an input (compact metric) space X of instances and output space $Y=\{-1,1\}$ of two labels aims at learning a (binary) classifier from samples that separate the instances in X into two classes. While the sampling process is under the control of a Borel probability measure $\rho$ on $Z:= X \times Y$, the performance of a classifier $C: X \rightarrow Y $ is assessed by the so-called misclassification error defined as the probability of the event $\{(x, y) \in X\times Y: C(x) \neq y\}$, that is
\begin{equation*}
    R(C):= \int_{X\times Y}I(-y, C(x))d\rho = Prob\{C(x)\neq y\}
\end{equation*}
where $I(a, b)= 1$ if $a = b$, and $0$ otherwise. The best classifier that minimizes the misclassification error is called a Bayes rule $f_c$ given by $f_c(x) = 1$ if $\rho(y=1|x) \geq  \rho(y= -1|x)$ and $-1$ otherwise, where $\rho(\cdot|x)$ denotes the conditional distribution of $\rho$ for $x \in X$. It turns out that for many convex loss functions $\phi: \mathbb{R} \rightarrow \mathbb{R^{+}}$, a Bayes rule can be expressed as $f_c = \mathrm{sgn}(f_{\rho}^{\phi})$ with $f_{\rho}^{\phi}$ being a minimizer of the generalization error
\begin{equation*}
    \mathcal{E}^{\phi}(f) =\int_{Z}\phi(yf(x))d\rho
\end{equation*}
over the set of measurable functions $f:  X \rightarrow \mathbb{R}$. Therefore, learning a Bayes rule is reduced to approximating $ f_{\rho}^{\phi}$ from hypothesis spaces. A special property of the classification
problem is that a Bayes rule $f_c$ taking binary values in ${1, -1}$ is often discontinuous, and $ f_{\rho}^{\phi}$ may also be discontinuous. For instance, for the hinge loss $\phi(v) = (1-v)_{+} := \max \{1-v, 0\}$, $f_{\rho}^{\phi}=f_c$. While this property motivates scholars to study approximation and learning in $L_p $ spaces with $1 \leq p \leq \infty $ of functions  $ f_{\rho}^{\phi}$ in Sobolev spaces $W^{r}_{p}$ with relatively small regularity index $r > 0$, this article, followed our novel approximation research of Korobov space \cite{Liu}, study the general error analysis in Korobov space.

The purpose of this paper is to \\
\\

Let $x=(x_1,\dots,x_d)\in\RR^d$ be the data vector and $m\in\NN$. A shallow neural network of width $m$ associated with a continuous activation function $\sigma: \RR \to \RR$ is defined by
\begin{equation}\label{defshallow}
f_m(x)=\sum\limits_{i=1}^m\beta_i\sigma(\alpha_i\cdot x-b_i),
\end{equation}
where $\{\alpha_i\}_{i=1}^m\subset\RR^d$ are connection vectors, $\{\beta_i\}_{i=1}^m\subset\RR$ weights and $\{b_i\}_{i=1}^m\subset\RR$ biases. 
The universality of shallow networks \cite{Cybenko,Hornik,Pinkus} asserts for any non-polynomial activation that any continuous function on any compact subset of $\RR^d$ can be approximated by output functions of the form (\ref{defshallow}) to an arbitrary accuracy when the number $m$ of hidden neurons is large enough. 
Rates of approximation by such output functions from the hypothesis space 
\begin{equation}\label{defspace}
H_m=\left\{\sum\limits_{i=1}^m\beta_i\sigma(\alpha_i\cdot x-b_i):\ \alpha_i\in\RR^d,\ \beta_i\in\RR,\ b_i\in\RR\right\}.
\end{equation}
were also studied in a large literature when $\sigma$ is replaced by a sigmoid type $C^\infty$ activation function.

\begin{definition}\label{defKorobov}
Let $D=[-1,1]^d$ with $d\in\NN$. The \textbf{Korobov space} $X^{2, p}(D)$ consists of functions $f\in\ C(D)$ vanishing on the boundary of $D$ and satisfying $D^{ k}f\in L_{p} (D)$ for any $k\in\ZZ_+^d$ with $\|k\|_\infty:=\max\{|k_1|,\dots, |k_d|\}\leq 2$.The norm is given by
\begin{equation}\label{KorobovNorm}
\lVert f\rVert_{X^{2, p}(D)}=\left\lVert\frac{\partial^{2d}f}{\partial x_1^2\dots\partial x_d^2}\right\rVert_{L_p (D)}+\lVert f\rVert_{L_p (D)}.
\end{equation}
\end{definition}

\section{ Main Results}
\subsection{Approximation rate of Korobov space using shallow ReLU networks}
In \cite{Liu}, we established a uniform approximation rate for approximating functions from the Korobov space $X^{2,\infty}([-1,1]^d)$ by ReLU shallow neural networks. In this paper, we extend the result to more general $L_p$ norm cases, that is, we proposed an upper bound of the rate of approximation from $X^{2, p}(D)$ by the same neural network structure.

\begin{theorem}\label{Thm1}
Let $d\in\NN$, $1\leq p \leq \infty$. Then there exists constants $C(d)$ and $C_5(d)$ depending only on $d$ such that for any $F\in X^{2,p}(D)$ and $m\in\mathbb N$, there holds
\begin{equation}\label{maininequation}
\begin{split}
\inf\limits_{f_m\in \mathcal{H}_m}\lVert F-f_m\rVert_{L_p(D)}\leq C(d)\|F\|_{X^{2,p}(D)}\sqrt{\log m}\left\{\begin{array}{ll}
        m^{-\frac{2(d+2)}{5d}}, &\quad 2\leq p\leq\infty,  \\
        m^{-\frac{2(d+2)}{5d+(\frac{2}{p}-1)d^2}}, & \quad 1\leq p<2,
    \end{array}\right.
\end{split}
\end{equation}
where
\begin{align}
    {\mathcal H}_{m}= \Big\{f_{m}: D \rightarrow \mathbb{R} |\ &f_{m}(x)= \sum_{k=1}^{m}\beta_{k}\sigma(\alpha_{k}\cdot x- b_{k}),  \ \|\alpha_{k}\|_{1} \leq 1,\nonumber \\
    & 0\leq b_{k} \leq 1, \ \|\beta_{k}\| \leq \frac{4\pi^2C_5(d)m^{\frac{1}{10}(1+\frac{2}{\alpha})}}{m}  \Big\}. \label{hypothesism}
\end{align}
\end{theorem}
A detailed proof will be given in the appendix.

\subsection{General error bound Classification }
The classification algorithm we study produces a regularized classifier
$\mathrm{sgn}(f_z)$ as the sign of $f_z$, a minimizer over the hypothesis
space ${\mathcal H}_{m}$ of the empirical error $\mathcal{E}_z(f)$ associated
with a convex loss function $\phi$ and a random sample $z :=\{(x_i, y_i)\}_{i=1}^N $ drawn according to $\rho$.

The function $f_z$ is the minimizer of regularization scheme
\begin{equation*}
f_z:= \arg \min_{f\in \mathcal{H}_m} \mathcal{E}_z( f ) = \arg \min_{f\in \mathcal{H}_m} \Big\{ \frac{1}{N} \sum_{i=1}^{N}\phi(y_if(x_i)) \Big\},
\end{equation*}
and our task is to estimate the excess misclassification error 
\begin{equation}\label{missclass}
    \mathcal{R} (\mathrm{sgn}(f_z)) -\mathcal{R}(f_c)
\end{equation}

Its convergence rates depend on the convexity of the loss $\phi$
and noise level of the underlying distribution $\rho$, which can be
measured simultaneously by the variancing power [33] of the pair $(\phi, \rho)$ defined as the maximum $ \tau \in [0, 1] $ such that for
some $C_1 > 0$
\begin{equation}\label{variance}
   \mathbb{E} \Big \{ (\phi(yf(x)) - \phi (yf^{\phi}_{\rho}(x)))^2 \Big \} \leq C_1 \Big \{ \mathcal{E}( f ) -\mathcal{E}(f_{\rho}^{\phi}) \Big \}^{\tau}
\end{equation}
holds for any measurable function $f : X \rightarrow R$.

We illustrate the learning rates for the $\eta$ norm loss $\phi(v)=(1 - v)_{+}^{\eta}$ with $\eta>1$ and the hinge loss with $\eta=1$.

For notation simplicity, the learning rates of the shallow neural network classification algorithm are given in terms of the approximation
error
\begin{equation*}
 D(H_m) = \inf_{f\in \mathcal{H}_m} \Big \{ \mathcal{E}( f ) -\mathcal{E}(f_{\rho}^{\phi}) \Big \}.
\end{equation*}
According to Theorem $4$ of \cite{D} and the fact that $\phi$ is $\eta-th$ time differentiable, the above equation yields

\begin{equation}
\inf_{f\in \mathcal{H}_m} \Big \{ \mathcal{E}( f ) -\mathcal{E}(f_{\rho}^{\phi}) \Big \} \leq  C^{\phi} ||f-f_{\rho}^{\phi})||^{\eta}_{L_{\rho X}^{\eta}},
\end{equation}
where $C^{\phi}=|\phi^{\eta}|_{L^{\infty}[-||f||_{\infty}-1, ||f||_{\infty}+1]}$.

\begin{theorem}\label{Thm2}
    Let $d\in\NN$, $1  \leq \eta \leq \infty$, $1  \leq p \leq \infty$, and $\phi(v)=(1 - v)_{+}^{\eta}$. If the pair $(\phi, \rho)$ has a variancing power $\tau \in [0,1]$ with (\ref{variance}) valid and the approximation error of the shallow hypothesis space satisfies
   \begin{equation}
\begin{split}
 D(H_m) =\inf_{f\in \mathcal{H}_m} \Big \{ \mathcal{E}( f ) -\mathcal{E}(f_{\rho}^{\phi})
 \Big \} \leq C^{\phi} \left\{\begin{array}{ll}
        m^{-\frac{2\eta(d+2)}{5d}}, &\quad 2\leq p\leq\infty,  \\
        m^{-\frac{2\eta(d+2)}{5d+(\frac{2}{p}-1)d^2}}, & \quad 1\leq p<2
    \end{array}\right.
\end{split}
\end{equation}
then by choosing 
\begin{equation}
  \begin{split}
   N= \left\{\begin{array}{ll}
        m^{(2-\tau)\frac{2\eta(d+2)+5d}{5d}}, &\quad 2\leq p\leq\infty\\
        m^{(2-\tau)\frac{2\eta(d+2)+5d+(2/p-1)d^2}{5d+(2/p-1)d^2}}, & \quad 1\leq p<2\\
    \end{array}\right.
  \end{split}  
\end{equation}

for any $\delta >0$, with confidence $1-\delta$, the excess misclassification error $\mathcal{R}(\mathrm{sgn}(f_z))-\mathcal{R}(f_c)$ of the induced classifier $\mathrm{sgn}(f_z)$ can be bounded as

\begin{equation}
  \begin{split}
   \left\{\begin{array}{ll}
       C_6 N^{-\frac{2\eta(d+2)}{(2-\tau)[2\eta(d+2)+5d]}}
       , &\quad 2\leq p\leq\infty, \eta>1\\
      C_6 N^{-\frac{2\eta(d+2)}{(2-\tau)(2\eta(d+2)+5d+(2/p-1)d^2)}}, & \quad 1\leq p<2,\eta>1\\
    \end{array}\right.
  \end{split}  
\end{equation}
where
$a=\max \{\eta2^{\eta-1}, C_1\}$
\begin{equation}
    C_6=\sqrt{[8C^{\phi}+24a+3C'_0+2(8C_1)^{\frac{1}{2-\tau}}][\max\{\log N, \log \frac{2N}{\delta}\}](\log\frac{2}{\delta})}
\end{equation}

and 

\begin{equation}
  \begin{split}
   \left\{\begin{array}{ll}
       C_7 N^{-\frac{2(d+2)}{(2-\tau)[2(d+2)+5d]}}
       , &\quad 2\leq p\leq\infty, \eta=1\\
      C_7 N^{-\frac{2(d+2)}{(2-\tau)(2(d+2)+5d+(2/p-1)d^2)}}, & \quad 1\leq p<2,\eta=1\\
    \end{array}\right.
  \end{split}  
\end{equation}
where
\begin{equation}
    C_7=[4C^{\phi}+24a+3C'_0+2(8C_1)^{\frac{1}{2-\tau}}][\max\{\log N, \log \frac{2N}{\delta}\}](\log\frac{2}{\delta}) 
\end{equation}

\end{theorem}

\subsection{Missclassification rate Under Noise Conditions}
In this subsection, we study the excess misclassification error under the Tsybakov noise condition.
\begin{definition}
     Let $\rho$ be a Borel probability measure. For the regression function $f_{\rho}: X\rightarrow \mathbb{R}$ is given by $f_\rho(x) = \eta(x)-(1-\eta(x))=2\eta(x)-1$ with $\eta(x): = \rho(y=1\mid x)$, 
     The Tsybakov noise condition with $\theta > 0$ for a Borel probability measure satisfies
   \begin{equation}\label{noise}
       \rho_{X}(\{x\in X \mid 0<|f_{\rho}(x)|\leq c_{\theta}m\}) \leq m^{\theta}, \quad \forall m > 0,
   \end{equation}
   where $c_{\theta}>0$ is a constant
\end{definition}
As a consequence, we define the Tsybakov function $T_{\rho}: [0,1] \rightarrow [0,1]$ equipped with  a Borel probability measure $\rho$ on $X \times Y$ as 
\begin{equation}
    T(r)=   \rho_{X}(\{x\in X \mid 0<|f_{\rho}(x)|\leq r\}), 
\end{equation}

\begin{theorem}\label{Thm3}
    Let $d\in\NN$, $\delta >0$ and $\phi(v)=(1 - v)_{+}^{2}$. If condition $(\ref{variance})$ valid with $\eta=2$ for some $r>0$ and Tsybakov noise condition $(\ref{noise})$ is satisfied, then by choosing 
\begin{equation}
  \begin{split}
   N= \left\{\begin{array}{ll}
        m^{\frac{9d+8}{5d}}, &\quad 2\leq p\leq\infty,\\
        m^{\frac{9d+8+(2/p-1)d^2}{5d+(2/p-1)d^2}}, & \quad 1\leq p<2,\\
    \end{array}\right.
  \end{split}  
\end{equation}
with confidence $1-\delta$, the excess misclassification error $\mathcal{R}(\mathrm{sgn}(f_z))-\mathcal{R}(f_c)$ of the induced classifier $\mathrm{sgn}(f_z)$ can be bounded as
\begin{equation}
  \begin{split}
   \mathcal{R}(\mathrm{sgn}(f_z))-\mathcal{R}(f_c)\leq C_{8} N^{-\frac{\theta}{2+\theta}}\left\{\begin{array}{ll}
       N^{\frac{4d+8}{9d+8}}
       , &\quad 2\leq p\leq\infty,\\
      N^{\frac{4d+8}{9d+8+(2/p-1)d^2}}, & \quad 1\leq p<2,\\
    \end{array}\right.
  \end{split}  
\end{equation}
where
\begin{equation}
    C_{8}=\sqrt{[8C^{\phi}+24\max \{4, C_1\}+3C_0+16C_1][\max\{\log N, \log \frac{2N}{\delta}\}](\log\frac{2}{\delta})}.
\end{equation}

\end{theorem}


\section{Proof of Main Result}
\subsection{General error bound}
In this subsection, we derive the learning rates stated in Theorem \ref{Thm2} for the 
excess misclassification error $\mathcal{R}(\mathrm{sgn}(f_z))-\mathcal{R}(f_c)$ of the induced classifier $\mathrm{sgn}(f_z)$. Since the stretch of the target function in the misclassification problem is $\{-1,1\}$, we introduce a truncated function of the minimizer function. We denote the truncate operator as
\begin{equation}
    \begin{split}
  \pi(f)(x)=
    \left\{\begin{array}{ll}
        1, &\quad \text{if } f(x)>1,  \\
       -1, & \quad  \text{if } f(x)<-1,\\
       f(x), & \quad \text{if } -1\leq f(x)\leq 1.
    \end{array}\right.
  \end{split}  
\end{equation}

\begin{theorem}\label{upperbound}
    Let $H$ be a compact subset of C(X) with $ B=\sup_{f\in H} ||f||_{\infty}$. Deﬁne $f_z$ by  with a convex loss function $\phi : \mathbb{R}\rightarrow \mathbb{R}^{+}$ and a random sample $z$. If $\phi(1)=0$ and the
pair $(\phi, \rho)$ has a variancing power $\tau \in [0,1]$ deﬁned by \ref{variance},
then for any $0 \leq \delta \leq 1$, with the probability at least $1-\delta$, the
excess generalization error $\mathcal{R}(\mathrm{sgn}(f_z))-\mathcal{R}(f_{\rho}^{\phi})$ can be bounded by
\begin{equation}
    4D(H_m) + \frac{8C'_0\log\frac{2}{\delta}}{3N} + 2\Bigg(\frac{8C_1\log\frac{2}{\delta}}{N}\Bigg)^{1/(2-\tau)} +24\epsilon^{*},
\end{equation}
where $C'_0 :=||\phi||_{L_\infty[-\max\{B,1\}, \max\{B,1\}]}$ and $\epsilon^{*}$ is the smallest positive number $\epsilon$ satisfying
$$\mathcal{N}\Bigg(H,\frac{\epsilon}{|\phi'_{+}(-1)|} \Bigg)\exp\Bigg\{-\frac{N\epsilon^{2-\tau}}{2C_1+\frac{4}{3}\phi(-1)\epsilon^{1-\tau}} \Bigg\} \leq \frac{\delta}{2}.$$

\end{theorem}

We need to bound the covering numbers to apply Theorem \ref{upperbound} to prove the main result stated in Theorem \ref{Thm2}. This step is given in the following lemma.

\begin{lemma}\label{coveringnumberlemma}
    Let ${\mathcal H}_{m}$ be the shallow neural network-generated hypothesis space defined in (\ref{hypothesism}). Then the covering number $\mathcal{N}(\epsilon,{\mathcal H}_{m} )$ satisfies
    $$ \log N(\epsilon,{\mathcal H}_{m}) \leq (d+2) m \log\frac{1}{\epsilon}+(d+2)m\Big(\log(1152e\pi^2C_5(d)) + \frac{d+2}{10d} \Big). 
    $$
\end{lemma}

In addition, by comparison theorem \cite{P}, \cite{D} and \cite{T}, we estimate the excess misclassification
error (\ref{missclass}) by the excess generalization error $ \mathcal{E}(\pi(f_z))- \mathcal{E}(f_{\rho}^{\phi})$
by taking $f = \pi (f_z)$.

\begin{lemma}\label{compare}
Let $f: X \rightarrow \mathbb{R}$ be a measurable function. For the $\eta$-norm loss $\phi(v) = (1-v)^{\eta}_{+}$ with $\eta>1$, there holds
\begin{equation}\label{c1}
     \mathcal{R} (\mathrm{sgn}(f_z)) -\mathcal{R}(f_c) \leq \sqrt{2(\mathcal{E}( f ) -\mathcal{E}(f_{\rho}^{\phi}))}
\end{equation}

For the hinge loss $\phi(v) = (1-v)_{+}$ we have  $f_{\rho}^{\phi}=f_{c}$ and

\begin{equation}\label{c2}
   \mathcal{R} (\mathrm{sgn}(f_z)) -\mathcal{R}(f_c) \leq \mathcal{E}( f ) -\mathcal{E}(f_{c}) 
\end{equation}
\end{lemma}

Now, we are in a position to prove Theorem \ref{Thm2}

\begin{proof}[Proof of Theorem \ref{Thm2}]
$ $\newline

By Lemma \ref{coveringnumberlemma} and $|\phi'_{+}(-1)|=\eta2^{\eta-1}$,
\begin{equation*}
    \log N(\frac{\epsilon}{|\phi'_{+}(-1)|},{\mathcal H}_{m}) \leq (d+2) m \log\frac{\eta2^{\eta-1}}{\epsilon}+(d+2)m\Big(\log(1152e\pi^2C_5(d)) + \frac{d+2}{10d} \Big).
\end{equation*}

Hence, for $0 \leq \delta \leq 1$, $\epsilon^{*}$ in equation can be bounded by solution 
$\tilde{\epsilon}$ of $h(\epsilon) \leq \log \frac{\delta}{2}$, where $h: \mathbb{R}^{+} \rightarrow \mathbb{R}$ is a decreasing function:

\begin{equation}
h(\epsilon)=(d+2)m\log(\frac{p2^{p-1}}{2}) +(d+2)m \log(1152e\pi^2C_5(d)+\frac{d+2}{10d})- \frac{N\epsilon^{2-\tau}}{2C_1+\frac{4}{3}\phi(-1)\epsilon^{1-\tau}}.
\end{equation} 

Let $A=d+2$, $Q=\log(1152e\pi^2C_5(d)+\frac{d+2}{10d})$, and note $|\phi'_{+}(-1)|=\eta2^{\eta-1}$, then
\begin{equation}
h(\epsilon)=Am\log(\frac{q2^{q-1}}{2}) + AQm- \frac{N\epsilon^{2-\tau}}{2C_1+\frac{2^{p+2}}{3}\epsilon^{1-\tau}}.
\end{equation}
Take $a=\max \{\eta2^{\eta-1}, C_1\}$, $\tilde{\epsilon}= a(\frac{\log \frac{2N}{\delta}}{N})^{\frac{1}{2-\tau}}m$, with the restriction
$$N \geq (\exp\{(d+2)m \log(1152e\pi^2C_5(d)\}(\frac{\delta}{2})^{2/(p-2)})^\frac{1}{2^{-p-2-(d+2)}},$$
we have
\begin{equation*}
h(\epsilon) \leq Am\log N+ AQm -\frac{1}{2^m}\log{2N}{\delta}m \leq \log \frac{\delta}{2}
\end{equation*}

We pick $\epsilon^{*}=\epsilon^{*}(m)= a(\frac{\log \frac{2N}{\delta}}{N})^{\frac{1}{2-\tau}}m $, by Theorem \ref{Thm1} and \ref{upperbound},

\begin{align*}
      &\mathcal{R}(\pi(f_z))-\mathcal{R}(f_{\rho}^{\phi}) \leq 
        4D(H_m) + \frac{8C'_0\log\frac{2}{\delta}}{3N} + 2(\frac{8C_1\log\frac{2}{\delta}}{N})^{\frac{1}{(2-\tau)}} +24\epsilon^{*}\\
        & \leq 4D(H_m) +\frac{3C_0}{N}\log\frac{2}{\delta} + 2(\frac{8C_1}{N})^{\frac{1}{(2-\tau)}}\log\frac{2}{\delta} +24a(\frac{\log \frac{2N}{\delta}}{N})^{\frac{1}{2-\tau}}m.\\
\end{align*}

If $2\leq p \leq \infty $, let $m=N^{\frac{1}{2-\tau}\frac{5d}{2q(d+2)+5d}}$
{\small\begin{equation*}
    \begin{split}
        &\mathcal{R}(\pi(f_z))-\mathcal{R}(f_{\rho}^{\phi})\\
        \leq& 4C^{\phi}N^{-\frac{2q(d+2)}{(2-\tau)[2q(d+2)+5d]}}+ 24a\left(\frac{\log(\frac{2N}{\delta})}{N}\right)^{\frac{1}{2-\tau}}N^{\frac{1}{2-\tau}\frac{5d}{2q(d+2)+5d}}\\
        &+\frac{3}{N}C'_0\log\frac{2}{\delta}+2\left(\frac{8C_1}{N}\right)^{\frac{1}{2-\tau}}\log\frac{2}{\delta}\\
        \leq&
    \left[4C^{\phi}+ 24a\left(\log\frac{2N}{\delta}\right)^{\frac{1}{2-\tau}}\right]N^{-\frac{2q(d+2)}{(2-\tau)[2q(d+2)+5d]}}+\left[3C'_0+2(8C_1)^{\frac{1}{2-\tau}}\right]\left(\log\frac{2}{\delta}\right)N^{-\frac{1}{2-\tau}}\\
    \leq& \left[4C^{\phi}+24a+3C'_0+2(8C_1)^{\frac{1}{2-\tau}}\right]\left[\max\left\{\log N,\log \frac{2N}{\delta}\right\}\right]\left(\log\frac{2}{\delta}\right)N^{-\frac{2q(d+2)}{(2-\tau)[2q(d+2)+5d]}}.
    \end{split}
\end{equation*}}
In the meanwhile,
If $1\leq p < 2 $, let $m=N^{\frac{1}{2-\tau}\frac{5d+(2/p-1)d^2}{2q(d+2)+5d+(2/p-1)d^2}}$,
{\small\begin{align*}
    &\mathcal{R}(\pi(f_z))-\mathcal{R}(f_{\rho}^{\phi}) \\ &
    \leq 4C^{\phi}N^{-\frac{2q(d+2)}{(2-\tau)[2q(d+2)+5d]}}+ 24a(\frac{\log(\frac{2N}{\delta})}{N})^{\frac{1}{2-\tau}}N^{\frac{1}{2-\tau}\frac{5d+(2/p-1)d^2}{2q(d+2)+5d+(2/p-1)d^2}}\\
    &+\frac{3}{N}C'_0(\log(\frac{2}{\delta}))+2(\frac{8C_1}{N})^{\frac{1}{2-\tau}}(\log(\frac{2}{\delta})))\\ &\leq
    [4C^{\phi}+ 24a(\log(\frac{2N}{\delta}))^{\frac{1}{2-\tau}}]N^{-\frac{2q(d+2)}{(2-\tau)(2q(d+2)+5d+(2/p-1)d^2)}}+[3C'_0+2(8C_1)^{\frac{1}{2-\tau}}](\log(\frac{2}{\delta}))N^{-\frac{1}{2-\tau}} \\&
    \leq [4C^{\phi}+24a+3C'_0+2(8C_1)^{\frac{1}{2-\tau}}][\max\{\log N, \log \frac{2N}{\delta}\}](\log\frac{2}{\delta})N^{-\frac{2q(d+2)}{(2-\tau)(2q(d+2)+5d+(2/p-1)d^2)}}.
\end{align*}}
together with (\ref{c1}) implies that for $\eta>1$, $2\leq p \leq \infty$, we have
\begin{equation}\label{pbig}
     \mathcal{R} (\mathrm{sgn}(f_z)) -\mathcal{R}(f_c) \leq C_6 N^{-\frac{2\eta(d+2)}{(2-\tau)[2\eta(d+2)+5d]}},
\end{equation}
and for $\eta>1$, $1\leq p <2$,
\begin{equation}\label{psmall}
     \mathcal{R} (\mathrm{sgn}(f_z)) -\mathcal{R}(f_c) \leq C_6N^{-\frac{2\eta(d+2)}{(2-\tau)(2\eta(d+2)+5d+(2/p-1)d^2)}},
\end{equation}
where
\begin{equation*}
    C_6= \sqrt{[8C^{\phi}+24a+3C'_0+2(8C_1)^{\frac{1}{2-\tau}}][\max\{\log N, \log \frac{2N}{\delta}\}](\log\frac{2}{\delta})}
\end{equation*}
\medskip
Similarly, together with (\ref{c2}) implies that for $\eta=1$, $2\leq p \leq \infty$
\begin{equation}
   \mathcal{R} (\mathrm{sgn}(f_z)) -\mathcal{R}(f_c) \leq   C_7 N^{-\frac{2(d+2)}{(2-\tau)[2(d+2)+5d]}}, 
\end{equation}
and for $\eta=1$, $1\leq p <2$,
\begin{equation}
   \mathcal{R} (\mathrm{sgn}(f_z)) -\mathcal{R}(f_c) \leq  C_7 N^{-\frac{2(d+2)}{(2-\tau)(2(d+2)+5d+(2/p-1)d^2)}},
\end{equation}
where
\begin{equation*}
 C_7= [4C^{\phi}+24\max \{1, C_1\}+3C'_0+2(8C_1)^{\frac{1}{2-\tau}}][\max\{\log N, \log \frac{2N}{\delta}\}](\log\frac{2}{\delta}).   
\end{equation*}
\end{proof}

\subsection{Improved rate Under Tsybakov Noise Conditions}
In this subsection, a proof of Theorem \ref{Thm3} is given, which is the ameliorated rate in Theorem \ref{Thm2} with Tsybakov noise conditions.
\begin{proof}[Proof of Theorem \ref{Thm3}]
According to Tsybakov condition(\ref{noise}), $T(l)= \mathcal{O})(l^{\theta})$. Let $l=\Gamma^{\theta/(2+\theta)}$ for some $\Gamma > 0$ yields
\begin{equation*}
    T(\sqrt{\Gamma/l})+l = \mathcal{O}(\Gamma^{\theta/(2+\theta)}).
\end{equation*}
Together with \cite[Proposition 1]{Huang} with $f=f_{z}$\\
\begin{equation} \label{enoise}
   \mathcal{R} (\mathrm{sgn}(f_z)) -\mathcal{R}(f_c) \leq \mathcal{O}(\mathcal{E}( \pi(fz) ) -\mathcal{E}(f_{\rho}^{\phi}))^{\frac{\theta}{2+\theta}}. 
\end{equation}
As illustrated in \cite{Q}, when $\eta=2$, i.e, the $2$-norm loss is considered, $\tau=1$ holds true. Combine (\ref{enoise}) and (\ref{pbig}) with $\eta=2 $ and $\tau=1$, we conclude
\begin{equation*}
    \mathcal{R} (\mathrm{sgn}(f_z)) -\mathcal{R}(f_c) \leq C_{8} N^{-\frac{\theta(4d+8)}{(2+\theta)(9d+8)}}
\end{equation*}
if $2\leq p \leq \infty$.

Similarly, for $1\leq p<2$, combine (\ref{enoise}) and (\ref{psmall}) with $q=2 $ and $\tau=1$,
\begin{equation*}
    \mathcal{R} (\mathrm{sgn}(f_z)) -\mathcal{R}(f_c) \leq    C_{8} N^{-\frac{\theta(4d+8)}{(2+\theta)[9d+8+(2/p-1)d^2)]}},
\end{equation*}
where
\begin{equation*}
    C_{8}= \sqrt{[8C^{\phi}+24\max \{4, C_1\}+3C'_0+16C_1][\max\{\log N, \log \frac{2N}{\delta}\}](\log\frac{2}{\delta})}.
\end{equation*}
\end{proof}

\section{Disucssion}

\begin{table}[h!]
\centering

\begin{tabular}{c c c c} 
\hline \hline
 Regularity & range & networks & error rate \\ [0.5ex] 
 \hline
 $f \in X^{2, p}([-1,1]^d)$  & $2 \leq r \leq 2d$ & deep ReLU & $N^{-2}$ \\ 
 \hline
$f \in X^{2, p}([-1,1]^d)$v& $r \leq \frac{d}{2}+2$ & DCNN & $N^{-2}$ \\
\hline
 $f \in X^{2, p}(D)$  & $2 \leq r \leq 2d$ & shallow ReLU & $N^{-\frac{2}{5}}, \quad 2\leq p\leq\infty$, \\ 
 
 $f \in X^{2, p}(D)$  & $2 \leq r \leq 2d$ & shallow ReLU &       $N^{-\frac{2}{5+(\frac{2}{p}-1)d}}, \quad 1\leq p<2$\\
\\ [1ex] 
 \hline
\end{tabular}
\caption{Approximation error in Korobov space}
\label{table:1}
\end{table}

\section{Appendix}
\begin{proof}[Proof of Theorem \ref{Thm1}]
We prove the Theorem in the appendix. \\
It has been shown in \cite[(3.8)]{MaoZhou} that
\begin{eqnarray}
    \left\|J_{N}(f)-f\right\|_{L_\infty (\TT^d)} &\leq& \sum_{j=1}^d 2^{2(j-1)}C_3\lVert f\rVert_{W^{2,\infty}(\TT^d)} N^{-2} \leq d4^{d}C_3\lVert f\rVert_{W^{2,\infty}(\TT^d)} N^{-2} \nonumber\\
    &\leq&d4^{d}C_1(d)C_3\|F\|_{X^{2,\infty}(D)}N^{-2}. \label{estimation1}
\end{eqnarray}


The following lemma is derived by taking $r=2$ in \cite[Lemma 1]{MaoZhouK}.
\begin{lemma}\label{BarronFourierseries}
For $k\in\ZZ^d$, let $\widehat{J_N}(k)$ be the Fourier coefficient of $J_{N}(f)$ at $k$ satisfying
$$J_{N}(f,x)=\sum\limits_{k\in\mathbb Z^d}\widehat{J_N}(k)e^{ik\cdot x}.$$
Then for each $m\in\NN$, there exists a function $f_m(x)=\sum\limits_{k=1}^m\beta_k\sigma(\alpha_k\cdot x-b_k)\in H_m$ such that
\begin{equation}\label{estimation2}
\lVert J_{N}(f)-f_m\rVert_{L_\infty(D)}\leq C_4v_{J_N,2}d^{3/2}\sqrt{\log m}m^{-\frac{1}{2}-\frac{1}{d}},
\end{equation}
where $C_4$ is an absolute constant,
\begin{equation}\label{defvf2}
v_{J_N,2}:=\sum\limits_{k\in\mathbb Z^d}\lvert\widehat{J_N}(k)\rvert\lVert k\rVert_1^2,
\end{equation}
and
$$\lvert\beta_k\rvert\leq\frac{4\pi^2v_{J_N,2}}{m},\qquad\lVert\alpha_k\rVert_1\leq1,\qquad 0\leq b_k\leq1,\qquad\forall1\leq k\leq m.$$
\end{lemma}

Let $L=\lceil\log_2N\rceil$, where $\lceil u\rceil$ denotes the smallest integer no less than $u>0$. For $\ell=(\ell_1,\dots,\ell_d)\in\{0,1,\dots,L\}^d$, let
$$\Lambda_\ell=\{k\in\ZZ^d:\ 2^{\ell_j-1}<|k_j|\leq2^{\ell_j},\ j=1,\dots,d\}.$$
By noticing the term vanishes at $k\in\ZZ^d$ with $k=0$ and $\|k\|_\infty>N$,
\begin{align}\label{upperbound v 1}
v_{J_N,2}=&\sum\limits_{\ell\in\{0,\dots,L\}^d}\sum\limits_{k\in\Lambda_\ell}\lvert  a_{k,2^L}\hat f(k)\rvert\lVert k\rVert_1^2\leq d\sum\limits_{\ell\in\{0,\dots,L\}^d}\sum\limits_{k\in\Lambda_\ell}\lvert  a_{k,2^L}\hat f(k)\rvert\left(k_1^2+\dots+k_d^2\right)\nonumber\\
=&d\sum\limits_{\ell\in\{0,\dots,L\}^d}\sum\limits_{k\in\Lambda_\ell}\sum\limits_{j=1}^d\lvert a_{k,2^L}\hat f(k)\rvert k_j^2.
\end{align}
For notation simplicity, we fix $\ell\in\{0,1,\dots,L\}^d$ and $j\in\{1,\dots,d\}$ and define
$$S_{\ell,j}:=\sum\limits_{k\in\Lambda_\ell}\lvert a_{k,2^L}\hat f(k)\rvert k_j^2=\sum\limits_{k\in\Lambda_\ell}\frac{\lvert a_{k,2^L}\hat f(k)\rvert\prod\limits_{s=1}^dk_s^2}{\prod\limits_{s\neq j}k_s^2}.$$
For $k\in\Lambda_\ell$, there holds $k_s^2\geq2^{\ell_s-1}$, hence
\begin{equation*}
\begin{split}
S_{\ell,j}\leq&\prod\limits_{s\neq j}2^{-2(\ell_s-1)}\left(\sum\limits_{k\in\Lambda_\ell}\lvert a_{k,2^L}\hat f(k)\rvert\prod\limits_{s=1}^dk_s^2\right)\\
\leq&\prod\limits_{s\neq j}2^{-2(\ell_s-1)}\sqrt{\left(\sum\limits_{k\in\Lambda_\ell}1^2\right)\left(\sum\limits_{k\in\Lambda_\ell}\left(\lvert a_{k,2^L}\hat f(k)\rvert\prod\limits_{s=1}^dk_s^2\right)^2\right)}\\
\leq&2^{2(d-1)}2^{\frac{\ell_j}{2}}\prod\limits_{s\neq j}2^{-\frac{3}{2}\ell_s}\sqrt{\sum\limits_{k\in\Lambda_\ell}\left(\lvert a_{k,2^L}\hat f(k)\rvert\prod\limits_{s=1}^dk_s^2\right)^2}.
\end{split}
\end{equation*}
Let
\begin{equation}\label{def TL}
T_{L}f(x)=\sum\limits_{\lVert k\rVert_\infty\leq2^L} \left( a_{k,2^L}\hat f(k)\prod\limits_{s=1}^dk_s^2\right)e^{ik\cdot x},\qquad x\in\TT^d,
\end{equation}
by Parseval's identity,
$$\sum\limits_{k\in\Lambda_\ell}\left(\lvert a_{k,2^L}\hat f(k)\rvert\prod\limits_{s=1}^dk_s^2\right)^2\leq\sum\limits_{k\in\ZZ^d}\lvert\widehat{T_{\ell}f}(k)\rvert^2=(2\pi)^{-d}\int_{\TT^d}\lvert T_L f(x)\rvert^2dx.$$
Hence
\begin{equation*}
\begin{split}
S_{\ell,j}\leq(2\pi)^{-\frac{d}{2}}2^{2(d-1)}2^{\frac{\ell_j}{2}}\prod\limits_{s\neq j}2^{-\frac{3}{2}\ell_s}\|T_Lf\|_{L_2(\TT^d)}.
\end{split}
\end{equation*}
Together with (\ref{upperbound v 1}),
\begin{align}\label{upperbound v 2}
v_{J_N,2}\leq&d\sum\limits_{\ell\in\{0,\dots,L\}^d}\sum\limits_{j=1}^dS_{\ell,j}\leq d2^{2(d-1)}(2\pi)^{-\frac{d}{2}}\|T_L f\|_{L_2(\TT^d)}\left(\sum\limits_{\ell\in\{0,\dots,L\}^d}\sum\limits_{j=1}^d2^{\frac{\ell_j}{2}}\prod\limits_{s\neq j}2^{-\frac{3}{2}\ell_s}\right)\nonumber\\
=&d2^{2(d-1)}(2\pi)^{-\frac{d}{2}}\|T_L f\|_{L_2(\TT^d)}\sum\limits_{j=1}^d\sum\limits_{\tau=0}^L2^{\frac{\tau}{2}}\sum\limits_{\substack{\ell\in\{0,\dots,L\}^d\\ \ell_j=\tau}}\prod\limits_{s\neq j}2^{-\frac{3}{2}\ell_s}\nonumber\\
\leq&d2^{2(d-1)}(2\pi)^{-\frac{d}{2}}\|T_L f\|_{L_2(\TT^d)}\sum\limits_{j=1}^d\sum\limits_{\tau=0}^L2^{\frac{\tau}{2}}\times2^{d-1}\nonumber\\
\leq&2d^22^{3(d-1)}(2\pi)^{-\frac{d}{2}}\|T_L f\|_{L_2(\TT^d)}2^{\frac{L}{2}}.
\end{align}
To prove an upper bound for $\|T_L f\|_{L_2(\TT^d)}$, write $\hat f(k)=(2\pi)^{-d}\int_{\TT^d}f(t)e^{-ik\cdot t}dt$ in the sum (\ref{def TL}), we have
\begin{align}\label{analyzesum2}
T_Lf(x)=&\sum\limits_{\|k\|_\infty\leq2^L} \left( a_{k,2^L}\hat f(k)\prod\limits_{s=1}^dk_s^2\right)e^{ik\cdot x}
\\
=&(2\pi)^{-d}\sum\limits_{\|k\|_\infty\leq2^L} a_{k,2^L}e^{ik\cdot x}\int_{\TT^d} \left(\prod\limits_{s=1}^dk_s^{2}\right) f(t)e^{-ik\cdot t}dt\nonumber\\
=&(2\pi)^{-d}\sum\limits_{\|k\|_\infty\leq2^L} a_{k,2^L}\int_{\TT^d}(-1)^{d}\frac{\partial^{2d}f}{\partial x_1^2\dots\partial x_d^2}(t)e^{ik\cdot(x-t)}dt\nonumber\\
=&(2\pi)^{-d}(-1)^{d}\int_{\TT^d}\frac{\partial^{2d}f}{\partial x_1^2\dots\partial x_d^2}(t)\prod\limits_{j=1}^d\left(\sum\limits_{k_j=-2^L}^{2^L}a_{k,2^L}^{[1]}e^{ik_j(x_j-t_j)}\right)dt\nonumber\\
=&(2\pi)^{-d}(-1)^{d}\int_{\TT^d}\frac{\partial^{2d}f}{\partial x_1^2\dots\partial x_d^2}(t)G_{2^L}(x-t)dt.
\end{align}
By Young's convolution inequality,
\begin{equation}
    \|T_L f\|_{L_2(\TT^d)}\leq(2\pi)^{-d}\left\{\begin{array}{ll}
        \left\|\frac{\partial^{2d}f}{\partial x_1^2\dots\partial x_d^2}\right\|_{L_{p}(\TT^d)}\lVert G_{2^L}\rVert_{L_{q}(\TT^d)}, & \quad 1\leq p<2, \\
        \left\|\frac{\partial^{2d}f}{\partial x_1^2\dots\partial x_d^2}\right\|_{L_{2}(\TT^d)}\lVert G_{2^L}\rVert_{L_1(\TT^d)}, &  \quad 2\leq p\leq\infty,
    \end{array}\right.
\end{equation}
where $q=\frac{2p}{3p-2}$ is the solution of $1+\frac{1}{2}=\frac{1}{p}+\frac{1}{q}$.

By definition (see also, e.g., \cite[(3.4)]{MaoZhou}), we have $\lVert G_{2^L}\rVert_{L_1(\TT^d)}\leq4^d$ and $\lVert G_{2^L}\rVert_{L_\infty(\TT^d)}\leq3C_2^d(2^{L+1}+1)^d$.\\

Then if $2\leq p\leq\infty$, we have
\begin{equation*}
    \begin{split}
    \|T_L f\|_{L_2(\TT^d)}\leq&(2\pi)^{-d}\left\|\frac{\partial^{2d}f}{\partial x_1^2\dots\partial x_d^2}\right\|_{L_{2}(\TT^d)}4^d.\\
    \end{split}
\end{equation*}
By H\"{o}lder inequality,\\
\begin{equation*}
\left\|\frac{\partial^{2d}f}{\partial x_1^2\dots\partial x_d^2}\right\|_{L_{2}(\TT^d)} \leq \left\|\frac{\partial^{2d}f}{\partial x_1^2\dots\partial x_d^2}\right\|_{L_{p}(\TT^d)} \left\|1\right\|_{L_{s}(\TT^d)},   
\end{equation*}
where $\frac{1}{s}=\frac{1}{2}-\frac{1}{p}$. Thus
\begin{equation*}
\begin{split}
  \|T_L f\|_{L_2(\TT^d)}\leq&4^d(2\pi)^{-d}\left\|\frac{\partial^{2d}f}{\partial x_1^2\dots\partial x_d^2}\right\|_{L_{p}(\TT^d)} \left\|1\right\|_{L_{s}(\TT^d)}\\
        \leq&4^d(2\pi)^{-d}\left\|\frac{\partial^{2d}f}{\partial x_1^2\dots\partial x_d^2}\right\|_{L_{p}(\TT^d)}(2\pi)^{d\left(\frac{1}{2}-\frac{1}{p}\right)}\\
       \leq&4^d(2\pi)^{(-\frac{1}{2}-\frac{1}{p})d}C_1(d)\|F\|_{X^{2,p}(D)}.
\end{split}
\end{equation*}
\medskip

If $1\leq p<2$, by Jensen's inequality,
$$\lVert G_{2^L}\rVert_{L_q(\TT^d)}\leq\lVert G_{2^L}\rVert_{L_1(\TT^d)}^{\frac{1}{q}}\lVert G_{2^L}\rVert_{L_\infty(\TT^d)}^{1-\frac{1}{q}}\leq3^{\frac{1}{2}}4^{d}C_2^{\frac{d}{2}}(2^{L+1}+1)^{\left(\frac{1}{p}-\frac{1}{2}\right)d}.$$
Thus
\begin{equation*}
    \begin{split}
        \|T_L f\|_{L_2(\TT^d)}\leq&(2\pi)^{-d}\left\|\frac{\partial^{2d}f}{\partial x_1^2\dots\partial x_d^2}\right\|_{L_{p}(\TT^d)}3^{\frac{1}{2}}4^{d}C_2^{\frac{d}{2}}(2^{L+1}+1)^{\left(\frac{1}{p}-\frac{1}{2}\right)d}\\
        \leq&3^{\frac{d+1}{2}}4^{d}C_2^{\frac{d}{2}}C_1(d)\|F\|_{X^{2,p}(D)}(2^{L+1}+1)^{\left(\frac{1}{p}-\frac{1}{2}\right)d}.
    \end{split}
\end{equation*}

\medskip
Combine above cases with (\ref{upperbound v 2}), we conclude
\begin{equation*}
    v_{J_N,2}\leq C_5(d)\|F\|_{X^{2,p}(D)}\left\{\begin{array}{ll}
        N^{\frac{1}{2}}, &\quad 2\leq p\leq\infty,  \\
        N^{\left(\frac{1}{p}-\frac{1}{2}\right)d+\frac{1}{2}}, & \quad 1\leq p<2.
    \end{array}\right.
\end{equation*}

where $C_5(d)=2\sqrt{3}d^22^{5d-3}\left(\frac{3C_2}{2\pi}\right)^{\frac{d}{2}}C_1(d)$.\\

Now Lemma \ref{BarronFourierseries} yields
\begin{equation}\label{estimation3}
\begin{split}
    &\inf\limits_{f_m\in H_m}\lVert J_{N}(f)-f_m\rVert_{L_\infty(D)}\\
    \leq&C_4C_5(d)d^{3/2}\|F\|_{X^{2,p}(D)}\sqrt{\log m}m^{-\frac{1}{2}-\frac{1}{d}}\left\{\begin{array}{ll}
        N^{\frac{1}{2}}, &\quad 2\leq p\leq\infty,  \\
        N^{\left(\frac{1}{p}-\frac{1}{2}\right)d+\frac{1}{2}}, & \quad 1\leq p<2.
    \end{array}\right.
\end{split}
\end{equation}
Together with (\ref{estimation1}), by taking $N=\lfloor m^{\frac{1}{5}(1+\frac{2}{d})}\rfloor$ for $2\leq p\leq\infty$ and $N=\left\lfloor m^{\frac{d+2}{5d+(\frac{2}{p}-1)d^2}}\right\rfloor$ for $1\leq p<2$, we conclude
\begin{equation}
\inf\limits_{f_m\in H_m}\lVert f-f_m\rVert_{L_p(D)}\leq C(d)\|F\|_{X^{2,p}(D)}\sqrt{\log m}\left\{\begin{array}{ll}
        m^{-\frac{2(d+2)}{5d}}, &\quad 2\leq p\leq\infty,  \\
        m^{-\frac{2(d+2)}{5d+(\frac{2}{p}-1)d^2}}, & \quad 1\leq p<2.
    \end{array}\right.
\end{equation}
where
$$C(d)=2d4^{d}C_1(d)C_3+d^{3/2}C_4C_5(d).$$
Thus we can complete the proof by noticing $f$ is an extension of $F$.
\end{proof}

\bibliographystyle{abbrvnat}

\end{document}